\documentstyle [11pt]{article}
\newfont{\bbb} {msbm10}

\begin{document}
\title{A Caveat on the Convergence of the Ricci Flow for Pinched Negatively Curved Manifolds}
\author{F. T. Farrell and P. Ontaneda\thanks{The first author was
partially supported by a NSF grant.
The second author was supported in part 
by a research  grant from CAPES, Brazil.}}
\date{November 5, 2003}
\maketitle

In his seminal paper \cite{H}, Hamilton initiated the Ricci flow method
for finding Einstein metrics on a closed smooth $n$-dimensional manifold $M^n$
starting with an arbitrary smooth Riemannian metric $h$ on $M^n$. He considered the evolution
equation $$ \frac{\partial}{\partial t}\, h=\frac{2}{n}\, r \, h \, -\, Ric$$

\noindent where $r =\int R\, d\mu /\,  \int d\mu$ is the average scalar curvature
($R$ is the scalar curvature) and $Ric$ is the Ricci curvature tensor of $h$.
Hamilton then spectacularly illustrated the success of this method by proving,
when $n=3$, that if the initial Riemannian metric has strictly positive
Ricci curvature it evolves through time to a positively curved Einstein
metric $h_\infty$ on $M^3$. And, because $n=3$, such a Riemannian metric
automatically has constant sectional curvature; hence $(M^3,h_\infty )$ is a spherical
space-form; i.e. its universal cover is the round sphere. Following Hamilton's approach
G. Huisken \cite{Hu}, C. Margerin \cite{Ma} and S. Nishikawa \cite{N} proved that, for every $n$,
sufficiently pinched to 1 $n$-manifolds (the pinching constant depending
only on the dimension)  can be deformed, through the Ricci flow, to a spherical-space form.\\

Ten years later R. Ye \cite{Y} studied the Ricci flow when the initial
Riemannian metric $h$ is negatively curved and proved that a negatively curved
Einstein metric is strongly stable; that is, the Ricci flow starting near such a
Riemannian metric $h$ converges (in the $C^\infty$ topology)
to a Riemannian metric isometric to $h$, up to scaling.
(We introduce the notation $h\equiv h'$ for two Riemannian metrics that are isometric
up to scaling.) 
In \cite{Y} R. Ye also proved that sufficiently pinched to -1 manifolds can be deformed,
through the Ricci flow, to hyperbolic manifolds, but the pinching constant in his theorem
depends on other quantities (e.g the diameter or the volume).
Ye's paper was motivated by the problem
on whether the Ricci flow can be used to deform every sufficiently pinched
to -1 Riemannian metric to an Einstein metric
(the pinching constant depending only on the dimension). 
His paper partially implements a scheme
proposed by Min-Oo \cite{M}.\\

In this short note we show that our previous results
\cite{FO} imply the existence of pinched negatively curved metrics for which the
Ricci flow does not converge smoothly. \\

{\it Definition.} We say the the Ricci flow for a negatively
curved Riemannian metric $h$ {\it converges smoothly}  if the Ricci flow, starting at $h$, 
is defined for all $t$ and converges (in the $C^\infty$ topology) to a well defined 
negatively curved (Einstein) metric.\\

It is a consequence of Ye's paper \cite{Y} that the Ricci flow for a negatively curved Riemannian metric 
$h$ converges smoothly if and only if the Ricci flow, starting at $h$,
eventually gets into a stable neighborhood of some negatively curved Einstein metric $h_{\infty}$.
By a stable neighborhood of $h_{\infty}$ we mean a neighborhood for which any Ricci flow
starting there converges to a metric isometric (up to scaling) to $h_\infty$. 
These stable neighborhoods can be taken as (sufficiently small) open sets just in the
$C^2$ topology.\\

{\bf Theorem.} {\it Given $n>10$ and $\epsilon >0$ there is a closed smooth
$n$-dimensional manifold $N$ such that

(i) $N$ admits a hyperbolic metric

(ii) $N$ admits a Riemannian metric $h$ with sectional curvatures in
$[-1-\epsilon ,-1+\epsilon]$ for which the Ricci flow does not converge smoothly.}\\

{\bf Proof.} Let $n>10$ and $\epsilon >0$. From \cite{FO} we have the following.\\

There are closed smooth manifolds $M_0$, $M_1$, $N$,
of dimension $n$, Riemannian metrics $g_0$, $g_1$ on $M_0$
and $M_1$, respectively, and smooth finite covers 
$p_0 : N\rightarrow M_0$,   $p_1 : N\rightarrow M_1$ such that: \\

(1) $M_0$ and $M_1$ are homeomorphic but not $PL$-homeomorphic.

(2) $g_0$ is hyperbolic

(3) $g_1$ has sectional curvatures in $[-1-\epsilon ,-1+\epsilon]$.

(4) There is a $C^\infty$ family of $C^\infty$ Riemannian metrics $h_s$ on $N$,
$0\leq s\leq 1$, with $h_0 =p_0 ^* g_0$, and $h_1 =p_1 ^* g_1$,
such that every $h_s$ has sectional curvatures in $[-1-\epsilon ,-1+\epsilon]$.\\

Note that $h_0$ is also hyperbolic. Now, since the Ricci flow preserves
isometries (see \cite{H}) we have that 
if the Ricci flow for $g_1$ does not converge smoothly,
then the Ricci flow for $h_1 =p_1 ^* g_1$ does not converge smoothly either, and we are done. Hence
we assume that the Ricci flow for $g_1$ converges smoothly. Let $g_{1,t}$ be the
Ricci flow starting at $g_{1,0}=g_1$, $0\leq t <\infty$, converging to the
negatively curved Einstein metric $g_{1,\infty }$. Note that, by Mostow' 
Rigidity Theorem and (1) above, $g_1$ and $g_{1,\infty}$
are non-hyperbolic. It follows that  $p_1^*g_1$ and $p_1^*g_{1,\infty}$
are also non-hyperbolic.\\

If the Ricci flow does not converge smoothly for some $h_s$, we are done. So, let us assume
that the Ricci flow converges smoothly for all $h_s$. We will show a contradiction. Write $h_{s,t}$,
for the Ricci flow starting at $h_{s,0}=h_s$, $0\leq t <\infty$, converging to the
negatively curved Einstein metric $h_{s,\infty }$. Then,
from the form of the evolution equation we have that $h_{1,t}=p_1^*g_{1,\alpha t}$,
for some constant $\alpha >0$, and
for all $0\leq t \leq \infty$. \\

{\bf Claim.} {\it  $(s,t)\mapsto h_{s,t}$ is continuous for $0\leq s\leq 1$, 
$0\leq t<\infty$, were we consider the space of Riemannian metrics with the
$C^\infty$ topology.}\\

To prove the claim we have to prove that the Ricci flow depends continuously
on the initial conditions. One way of doing this directly is by using Hamilton's
proof of the local-in-time existence and uniqueness of the Ricci flow
(see \cite{H}). Let $f_0$, $\bar f$, $\bar h$ be as in the 
the proof of theorem 5.1 of \cite{H}, p.263. Let $f'_0$ be another initial condition. 
If $f'_0$ is close to $f_0$ (in the $C^\infty$ topology)
then we can find a $\bar f'$ close to $\bar f$ ($\bar f'$ with the same properties as
$\bar f$, but with respect to $f'_0$). Then $\bar h'$ is close to $\bar h$, where
$\bar h'$ is defined in a similar way as $\bar h$. Since the inverse function is
continuous, follows that $f'$ and $f$ are close, where  $f'$ and $f$ are the inverses
of some $h'$ and $h$ (which are chosen close to $\bar h'$ and $\bar h$ and vanishing
on some small interval $[0,\epsilon )$. This proves the claim.\\

Since $h_{0}$ is hyperbolic we have that $h_{0,t}=h_0$
for all $0\leq t\leq \infty$. Since every negatively curved Einstein metric 
is stable (see \cite{B}, p.357)
we can assume that all negatively curved Einstein metrics $h_{s,\infty}$
have neighborhoods $V_{s}$ for which any Ricci flow starting in $V_s$
converges to a metric isometric (up to scaling) to $h_{s,\infty}$ (see \cite{Y}, p.873) 
(in particular, $h_{s,\infty}\equiv h_0$, for sufficiently small $s$). 
It follows that every $s\in [0,1]$ has an open neighborhood $I_s$ such that
$h_{s',\infty}\equiv h_{s,\infty}$ for all $s'\in I_s$.
Then the map $s\mapsto [h_{s,\infty}]$ from [0,1] to ${\cal M}_N / \equiv$
is locally constant and hence continuous. (Here $[h_{s,\infty}]$ denotes the equivalence
class of $h_{s,\infty}$ in the quotient space ${\cal M}_N / \equiv$
of all isometry classes of Riemannian metrics on $N$.)
This is a contradiction because $h_{0,\infty}=h_0$ is hyperbolic and $h_{1,\infty}$ 
is not hyperbolic. This proves the theorem.\\
\vspace{.3in}

Recall that ${\cal M}_P$ denotes the space of all Riemannian metrics on a smooth manifold $P$.
For $\epsilon >0$, let ${\cal M}^{\epsilon}_P$ denote the space of $\epsilon$-pinched to -1
Riemannian metrics on $P$. Also, ${\cal E}_P\subset {\cal M}_P$ will denote the space of negatively curved 
Einstein metrics on $P$. Recall that ${\cal E}_P/\equiv$ is discrete, see \cite{B}, p.357.\\

{\it Definition.} Let $\epsilon > 0$ and $n$ be a positive integer. An {\it Einstein correspondence} 
$\Phi :{\cal M}^\epsilon\rightarrow {\cal E}$ {\it for $n$-dimensional manifolds}  
is a family of maps $\Phi_P :{\cal M}^\epsilon_P\rightarrow{\cal E}_P$, for every $n$-dimensional 
manifold $P$ for which  ${\cal M}^\epsilon_P$ is not empty.
We say that $\Phi$ is {\it cover-invariant} if $\Phi ( p^*g)\, =\, p^*(\Phi (g))$ for every finite cover
$p: P\rightarrow Q$ and $g\in{\cal M}^\epsilon_Q$, for which $\Phi_Q$ is defined.\\

We say that $\Phi$ is {\it continuous} if each $\Phi_P :{\cal M}^\epsilon_P\rightarrow{\cal E}_P$
is continuous. Here we consider ${\cal M}^\epsilon_P$ with the $C^\infty$ topology and ${\cal E}_P$
with the $C^2$ topology.\\

Let $h$, $h'$ $\in {\cal M}_P$. Write $h\equiv_0h'$ provided $(P,h)$ is isometric to $(P,h')$, up to scaling,
via an isometry homotopic to $id_P$. Notice that the fibers of ${\cal E}_P/\equiv_0\,\, \rightarrow \,\, {\cal E}_P/\equiv$
are discrete; and hence ${\cal E}_P/\equiv_0$ is also discrete.\\

The following corollary is a direct consequence of the theorem above.\\

{\bf Corollary 1.} {\it 
Suppose that there are  $\epsilon >0$ and $n>10$ for which there exists a cover-invariant
Einstein correspondence $\Phi$. Then there is a closed
$n$-dimensional Riemannian manifold $N$, with metric $h\in{\cal M}^\epsilon_N$, for which 
the Einstein metric $\Phi ( h)$ is unreachable by the Ricci flow starting at $h$.}\\

{\bf Corollary 2.} {\it 
Suppose that there are  $\epsilon >0$ and $n\geq 6$ for which there exists an
Einstein correspondence $\Phi$. Then there is a closed
$n$-dimensional manifold $N$ that admits, at least, two non-isometric
negatively curved Einstein metrics. Moreover, one metric can be chosen to be 
hyperbolic.}\\

{\bf Proof.} 
From \cite{FO2} we have the following.\\

There are closed connected smooth manifolds $M_0$, $M_1$, $N$,
of dimension $n$, Riemannian metrics $g_0$, $g_1$ on $M_0$
and $M_1$, respectively, and smooth two-sheeted covers 
$p_0 : N\rightarrow M_0$,   $p_1 : N\rightarrow M_1$ such that: \\

(1) $M_0$ and $M_1$ are homeomorphic but not $PL$-homeomorphic.

(2) $g_0$ is hyperbolic

(3) $g_1$ has sectional curvatures in $[-1-\epsilon ,-1+\epsilon]$.\\

Then the two non-isometric negatively curved Einstein metrics on $N$ are
$p_0^*(g_0)$ and $p_1^* (\Phi (g_{1}) )$. This proves the corollary.\\

The more general form of the following corollary was suggested to us by Rugang Ye.\\

{\bf Corollary 3.} {\it A cover-invariant Einstein correspondence cannot be
continuous.}\\

{\bf Remark.} Note that we are not assuming that $\Phi$ fixes hyperbolic metrics.
If we assumed that $\Phi (hyperbolic\,\,\, metric)=(hyperbolic\,\,\, metric),$ the proof
of the corollary would be much easier.\\

{\bf Proof.} We use the notation from the proof of the theorem.
Let us suppose that there exists a continuous cover-invariant Einstein correspondence.
We will show a contradiction.\\

Let $G_i\subset Diff\, (N)$, be (finite) subgroups of the group $Diff\, (N)$, of all self-diffeomorphisms
of $N$, such that $N/G_i=M_i$, $i=0,1$. 
Note that $G_i\subset Iso(N,\Phi(h_i))$, since $\Phi (h_i)=p_i^*(\Phi(g_i))$, 
where $Iso\, (N,\Phi (h_i))\subset Diff\, (N)$ is the subgroup consisting 
of all isometries of the negatively curved Einstein manifold $(N,\Phi (h_i))$.
Let $Top\, (N)$ and $Out\, (\pi_1 N)$ denote
the group of all self-homeomorphisms of $N$ and the group of outer automorphisms of $\pi_1(N)$, 
respectively.
Recall that $Out(\pi_1N)$ can be identified with $\pi_0(E(N))$, where $E(N)$ is the $H$-space 
consisting of all self-homotopy equivalences of $N$. (This is because $N$ is aspherical.)
We have the following diagram of group homomorphisms
$$Diff \, (N) \stackrel{\alpha}{\rightarrow}Top\, (N)\stackrel{\beta}{\rightarrow}Out\, (\pi_1 N)$$
where $\alpha$ is the inclusion, and $\beta$ is the composition of the inclusion
$Top\, (N)\rightarrow E(N)$
and the quotient map $E(N)\rightarrow \pi_0(E(N))$. Write $\gamma=\beta\alpha$.

It was shown in \cite{FO2}, \cite{FO}, that  $G_0$ and $G_1$ are conjugate in $Top\, N$, 
via a homeomorphism homotopic to $id_N$; hence $\gamma G_0=\gamma G_1$.\\

Now, since $\Phi$ continuous and ${\cal E}_N/\equiv_0$ is discrete, the composition 
$$[0,1]\,\, \stackrel{h_t}{\rightarrow} \,\, {\cal M}_N^\epsilon \,\, \stackrel{\Phi}{\rightarrow} 
\,\, {\cal E}_N\stackrel{}{\rightarrow}\,\,  {\cal E}_N/\equiv_0$$

\noindent must be constant; hence $\Phi (h_0)\equiv_0\Phi(h_1)$.
It follows that  $G_1$ is conjugate in $Diff\, (N)$
to a subgroup of $Iso\, (N,\Phi (h_0))$ 
via a diffeomorphism $f$ homotopic to $id_N$; i.e. $f^{-1}G_1f\subset Iso\, (N,\Phi (h_0))$. 
Note that  $\gamma (f^{-1}G_1 f)=\gamma (G_1)$ since $f\sim id_N$; hence 
$\gamma (f^{-1}G_1 f)=\gamma (G_0)$. This implies that  $f^{-1}G_1 f=G_0$ since both
$f^{-1}G_1 f$ and $G_0$ are subgroups of $Iso\, (N,\Phi (h_0))$ and Borel-Conner-Raymond
showed (see \cite{CR}, p.43)
that $\gamma$ restricted to compact subgroups of $Diff\, (N)$ is monic. (Recall that $N$ is aspherical
and the center of $\pi_1(N)$ is trivial.)
It follows that $f$ induces a diffeomorphism between $M_0=N/G_0$ and 
$M_1=N/G_1$, which is a contradiction. This proves the corollary.\\

{\it Definition.} We say the the Ricci flow for a  
negatively curved Riemannian metric $h$ {\it converges smoothly (weakly) } if the Ricci flow, starting at $h$, is
defined for all $t$ and converges (in the $C^\infty$ topology) to a well defined (Einstein) metric
(not necessarily negatively curved).\\

The following is also a corollary of the proof of the theorem above.\\

{\bf Corollary 4.} {\it  Given $n>10$ and $\epsilon >0$ there is a closed
$n$-dimensional manifold $N$ that admits a hyperbolic metric $h_0$
and a Riemannian metric $h_a$, with sectional curvatures in
$[-1-\epsilon ,-1+\epsilon]$, that satisfies the following. Either the Ricci flow for $h_a$
does not converge (weakly) or $N$ supports a  non-stable (hence not negatively curved)
Einstein metric $\tilde h$ satisfying:

\noindent(i)  There is a $C^\infty$ family of $C^\infty$ Riemannian metrics $h_s$ on $N$,
$0\leq s\leq a$, such that every $h_s$ has sectional curvatures 
in $[-1-\epsilon ,-1+\epsilon]$ and:

\hspace{.1in} (a) the Ricci flow, starting at $h_s$, converges to a  metric isometric 

\hspace{.1in} to the hyperbolic metric $h_0$ provided 
$0\leq s<a$.

\hspace{.1in} (b) the Ricci flow, starting at $h_a$, converges smoothly (weakly)  to

\hspace{.1in} the non-stable not negatively curved Einstein metric $\tilde h$.

\noindent(ii) There is a sequence of metrics $h_n$ converging (in the $C^\infty$topology)
to $\tilde h$ such the the Ricci flow, starting at each $h_n$, converges to a metric 
isometric to the hyperbolic metric $h_0$.}\\

{\bf Proof.} 
We use all notation from the proof of the theorem. 
As before, if the Ricci flow for $g_1$ does not converge smoothly (weakly) we are done. Let us assume then that
the Ricci flow for $g_1$ and all $h_s$ converges smoothly (weakly).
Let 

$$a\, =\, sup\, \{\, s\in[0,1]\,:\,\, h_{s',\infty}=h_{0},\,\, s'\in [0,s] \,\}$$

Then $0<a\leq 1$. It follows that $h_{a,\infty}$ is a  non-stable
(hence not negatively curved) Einstein metric. Take ${\tilde h}=h_{a,\infty}$. This proves part 
{\it (i)}.\\

To prove {\it (ii)}, note that we can choose a sequence $s_n$ with $s_n\rightarrow a$, $s_n <a$,
such that $h_n=h_{s_{n},n}\rightarrow \tilde h$.
Then the Ricci flow, starting at $h_n$, converges to $h_{s_{n},\infty}\equiv h_0$. This proves the corollary.\\

{\bf Acknowledgments.} We are very grateful to Rugang Ye and Michael Anderson for their helpful comments.
We wish to thank C.S. Aravinda for the useful information he provided to us.

F.T. Farrell

SUNY, Binghamton, N.Y., 13902, U.S.A.\\

P. Ontaneda

UFPE, Recife, PE 50670-901, Brazil

\end{document}